\documentclass[conference]{IEEEtran}
\usepackage{cite}

\ifCLASSINFOpdf
  \usepackage[pdftex]{graphicx}
\else
\fi
\usepackage{bm}

\usepackage{amsmath}
\usepackage{amsfonts}

\usepackage{array}
\usepackage{url}

\begin{document}
%
\title{A Julia Module for Polynomial Optimization with Complex Variables applied to Optimal Power Flow}


\author{\IEEEauthorblockN{Julie Sliwak\IEEEauthorrefmark{1}\IEEEauthorrefmark{2}\IEEEauthorrefmark{4},
Manuel Ruiz\IEEEauthorrefmark{1},
Miguel F. Anjos\IEEEauthorrefmark{3}, 
Lucas L\'etocart\IEEEauthorrefmark{4} and
Emiliano Traversi\IEEEauthorrefmark{4}}
\IEEEauthorblockA{\IEEEauthorrefmark{1}RTE R$\&$D, Paris La D\'efense, France}
\IEEEauthorblockA{\IEEEauthorrefmark{2}D\'epartement de math\'ematiques et de g\'enie industriel, Polytechnique Montr\'eal, Montr\'eal (Qu\'ebec), Canada}
\IEEEauthorblockA{\IEEEauthorrefmark{3}School of Mathematics, University of Edinburgh, Edinburgh, Scotland, UK EH9 3FD}
\IEEEauthorblockA{\IEEEauthorrefmark{4}LIPN, UMR CNRS 7030, Universit\'e Paris 13, Sorbonne Paris Cit\'e, Villetaneuse, France}}



\maketitle

\begin{abstract}
Many optimization problems in power transmission networks can be formulated as polynomial problems with complex variables. A polynomial optimization problem with complex variables consists in optimizing a real-valued polynomial whose variables and coefficients are complex numbers subject to some complex polynomial equality or inequality constraints. These problems are usually directly converted to real variables, either using the polar form or the rectangular form. In this work, we propose a Julia module allowing the representation of polynomial problems in their original complex formulation. This module is applied to power systems optimization and its generic design enables the description of several variants of power system problems. Results for the Optimal Power Flow in Alternating Current problem and for the Preventive-Security Constrained Optimal Power Flow problem are presented.
\end{abstract}


%
\IEEEpeerreviewmaketitle

\section{Introduction}

Transmission System Operators (TSOs) like France's RTE must take into account resistive losses in their models because of the Joule effect. The Direct Current (DC) description defines linear models that often neglect the resistive losses, while the Alternating Current (AC) description enables the representation of resistive losses plus the definition of a voltage plan. However, AC models are much more difficult to solve.

AC problems in power networks are nonlinear optimization problems with complex variables as AC modelling involves complex quantities like the voltage or the current. More specifically, they are part of the class of Polynomial Optimization Problems with Complex Variables $(POP-\mathbb{C})$, highlighted in the recent works of Josz and Molzahn \cite{josz2018lasserre} on the application of the Lasserre hierarchy to $(POP-\mathbb{C})$. 

To the best of our knowledge, there is at present no tool to represent $(POP-\mathbb{C})$ problems in their original formulation and they are usually directly converted into problems with real variables using polar or rectangular form. Even though some tools, like the well-known power system package MATPOWER \cite{matpower}, exploit the complex structure in the autodifferentiation procedure, the problems are defined with real variables. A software for $(POP-\mathbb{C})$ would allow to represent problems in a generic form and to test and develop specific methods for problems with complex variables (Lasserre hierarchy \cite{josz2018lasserre}, solver for SemiDefinite Programming problems with complex variables \cite{jcg_SDP_solver_complex_numbers}).

Another challenge for power system modelling is the need for modelling flexibility since optimization models become more complicated within the context of energy transition. Power system software have been designed for conventional networks and have now to adapt to intermittent ressources and emerging technologies like storage. Indeed, most power systems software use predefined structures, which has the advantage of being easier to use but does not allow users to model new problems. The recent Julia module PowerModels.jl \cite{coffrin2018powermodels} takes a step in the direction of modelling flexibility by providing several formulations of power flow problems but all models are predefined and expressed with real variables. The Python tool pandapower\cite{pandapower2018} allows even more flexibility by providing element-based models. Although there are many possible elements, they are predefined and their attributes are expressed with real numbers.

The Grid Optimization Competition launched by ARPA-E (Advanced Research Projects Agency-Energy) in 2017 is a first step towards network transformation \cite{GOC}. The main objective of this competition is to accelerate the development of new methods for solving power system problems in modern grids. In this competition, modelling flexibility is required both for the definition of network elements and for the mathematical formulation. For example, the problems involve logical constraints. One can use binary variables, if supported, to model this type of constraints.

In this work, we propose an open-source Julia \cite{Julia} module for Mixed-Integer $(POP-\mathbb{C})$ allowing the representation of power system problems in their original formulation \cite{MathProgComplex.jl}. From this module, we have designed a generic power system tool which enables the user to define the networks elements directly. We demonstrate the convenience of this module by applying it to two important problems in power systems: the Optimal Power Flow problem in Alternating Current (ACOPF)\cite{carpentier,HistoryofOPF,OPF_biblioI, OPF_biblioII} and the Preventive Security Constrained Optimal Power Flow (PSCOPF) problem \cite{PSCOPF}. 

The paper is organized as follows. Section II defines Mixed-Integer $(POP-\mathbb{C})$ and presents the Julia module we propose. Section III demonstrates the application of our tool to ACOPF and PSCOPF problems. Future research directions are exposed in section IV and section V concludes the paper.

\section{Julia Module MathProgComplex.jl}

\subsection{Motivation and background}

A polynomial optimization problem with complex variables $(POP-\mathbb{C})$ consists in optimizing a real-valued multivariate complex polynomial subject to some real-valued complex polynomial equality or inequality constraints. A real-valued multivariate complex polynomial is a polynomial whose variables and coefficients are complex numbers but whose value is always real. These problems can be extended to Mixed-Integer problems $(MIPOP-\mathbb{C})$ in which some of the variables are real and integer. We express $(MIPOP-\mathbb{C})$ as:

\begin{equation}\tag{$MIPOP-\mathbb{C}$}
\begin{aligned}
& \underset{z}{\text{min}}
& & \sum_{\alpha,\beta} p^0_{\alpha,\beta} z^{\alpha} \overline{z}^{\beta} \\
& \text{s.t.}
& & \sum_{\alpha,\beta} p^i_{\alpha,\beta} z^{\alpha} \overline{z}^{\beta} \geq 0, \; i = 1, \ldots, m.\\
& & &z \in \mathbb{C}^n, z_k\in \mathbb{N}\ \forall k \in K\subset\{1,...,n\}
\end{aligned}
\end{equation} where $\overline{z}$ represents the conjugate of $z$, the sums over $\alpha,\beta\in \mathbb{N}^n$ are finite and the coefficients $p^i_{\alpha,\beta} \in \mathbb{C}^1$.

There are currently very few methods to solve a $(MIPOP-\mathbb{C})$ problem in complex variables directly; the problem has to be converted into a problem with real variables to be solved. There are two possibilities for each variable: either using the polar form which implies nonlinear expressions with trigonometric functions or using the rectangular form which leads to polynomial expressions. In either case the result is a mixed-integer nonconvex problem.  

There are two categories of methods for mixed-integer nonconvex optimization: heuristics and exact methods. The goal of heuristics is to find good feasible solutions but without guarantee of optimality. In practice, the best known heuristic for power system problems, such as ACOPF problems, is a nonlinear interior point method. Still, other methods like Sequential Quadratic Programming or nondeterministic methods (genetic algorithms for example) can also be applied \cite{OPF_biblioI,OPFsurvey_deterministic_and_stochastic_methods}. On the other hand, exact methods aim at proving global optimality. Two of the most common exact algorithms for mixed integer problems are spatial Branch-and-Bound and Branch-and-Cut algorithms. Both algorithms consist in splitting the feasible domain and computing lower bounds at every node by solving convex relaxations \cite{B&C_for_nonconvex_QCQP_varCbounded,SOCP_relax_Sun_2017,castro2016spatial}. Moreover, it is also possible to use convergent convexification techniques \cite{kuang2017alternative} such as the Lasserre hierarchy \cite{molzahn2014moment,josz_SOS_OPF,Lasserre_partial_sparsity}, consisting in computing tighter and tighter conic relaxations until global optimality of a solution is proven. 

Our module MathProgComplex.jl provides a structure and methods to work with $(MIPOP-\mathbb{C})$ so that problems can be represented in their original formulation with complex variables. For now, only binary variables and not integer ones are supported. Basic operations such as addition, multiplication, conjugation or modulus are implemented so that problems can be expressed in a form easier to understand. An additional advantage is that there is no need to choose a real representation from the beginning and new methods for complex problems can be tested. Methods for real problems can also be applied by converting complex problems to real problems and this conversion can be done at the last moment, i.e., when problems have to be solved. This late conversion to real variables will enable a more convenient comparison between different real formulations of a common complex problem. For now, our module provides a function to convert problems using rectangular form but other conversion functions will be implemented. To give an overview of what can be done using MathProgComplex.jl, some examples are given in the following section.

\subsection{Examples}

The base structure is \textit{Variable}, which is a pair (string, type). A variable can be of \textit{Complex, Real or Bool} type. The string simply defines the name of the variable. The \textit{Exponent} and \textit{Polynomial} structures are constructed from \textit{Variables} by calling the respective constructors or with algebraic operations. There is also a \textit{Point} structure holding the variables at which polynomials can be evaluated.

For example, to define the complex polynomial ~$(1+4im)*x^2\overline{y}^3+3xy+bx$ with $b$ a binary variable, we can proceed as follows:

\begin{equation*}
\begin{array}{l}
x=Variable("x", Complex)\\
y=Variable("y", Complex)\\
b=Variable("b", Bool)\\
expo1=x^2*conj(y)^3\\
expo2=x*y\\
expo3=b*x\\
p=(1+4im)*expo1+3*expo2+expo3\\
\end{array}
\end{equation*} where $conj(z)$ stands for the conjugate of $z\in \mathbb{C}$

or more succinctly:
\begin{equation*}
\begin{array}{l}
x=Variable("x", Complex)\\
y=Variable("y", Complex)\\
b=Variable("b", Bool)\\
p=(1+4*im)*x^2*conj(y)^3+3*x*y+b*x\\
\end{array}
\end{equation*}

It is also easy to evaluate this polynomial at point $(1+2im, 3im, 0)$:

\begin{equation*}
\begin{array}{l}
pt=Point([x,y,b],[1+2*im,3*im,0])\\
evaluate(p,pt)\\
>-145+28im
\end{array}
\end{equation*}

As a constraint is simply defined by a polynomial and complex bounds, polynomial problems are straightforward to implement. For instance, the problem \eqref{P}: 

\begin{equation}\label{P}\tag{P}
\begin{aligned}
& \underset{x}{\text{min}}
& & \frac{1}{2}(x+\overline{x}-im (x-\overline{x})) \\
& \text{s.t.}
& & |x|^2=1\\
& & &x \in \mathbb{C}
\end{aligned}
\end{equation}
is defined as:
\begin{equation*}\label{codeP}
\begin{array}{l}
p=Problem()\\
x=Variable("x", Complex)\\
set\_objective!(pb,0.5*(x+conj(x)-im*(x-conj(x)))\\
add\_constraint!(pb,abs2(x)==1)
\end{array}
\end{equation*} where $abs2(z)$ stands for the squared modulus of $z\in \mathbb{C}$

To solve a complex problem like \eqref{P}, one has to convert it into a real problem. Our module provides the function $pb\_cplx2real(\cdot)$ to convert problems with complex variables into problems with real variables using the rectangular form. For instance, ($pb\_cplx2real(p)$) returns the following problem:

\begin{equation}\label{Preal}\tag{Preal}
\begin{aligned}
& \underset{x_{Re}, x_{Im}}{\text{min}}
& & x_{Re}+x_{Im} \\
& \text{s.t.}
& & x_{Re}^2+x_{Im}^2=1\\
& & &x_{Re},x_{Im} \in \mathbb{R}
\end{aligned}
\end{equation} Then the problem can either be converted into a JuMP \cite{JuMP} model or be exported into a dedicated text format from which, for example, an AMPL \cite{AMPL} model can be created.

In the next section, we present two power system problems that can be tackled using MathProgComplex.jl.

\section{Application to AC Optimal Power Flow}

\subsection{AC Optimal Power Flow}
The objective of the ACOPF problem is to provide an optimal dispatch meeting the electrical demand, satisfying physical laws and engineering constraints while minimizing the generation costs. Several formulations of the ACOPF problem can be found in the literature \cite{HistoryofOPF,IVformulation_LinearApproximation}. Some of them differ by the constraints included: thermal limits on lines can either be neglected or be modeled in terms of power or current. Others differ by the representation of the complex variables, which means that they are different from the real point of view but not from the complex point of view. Indeed, for a given ACOPF complex problem, several real problems can be defined depending if the polar or the rectangular form is used for each type of variable. Our module allows the representation of any type of ACOPF problem with complex variables.

For clarity, let us define a typical formulation of an ACOPF problem using our module. First, power transmission networks can be modelled as oriented graphs $T=(N,B)$ in which buses represent demand and/or generation points and where arcs represent transmission lines or transformers. We express a typical ACOPF problem as following:

\begin{equation}\label{ACOPF}\tag{ACOPF}
\begin{aligned}
& \underset{V, S^{gen}}{\text{min}}
& & \sum_{g\in G} c_g(Re(S_g^{gen}))\\
& \text{s.t.}
& & S^{gen}_n = S^{load}_n+\mathbf{S^{shunt}_n} +\\
& & &\sum_{b\in B^-(n)}S^{d}_b(V)+\sum_{b\in B^+(n)}S^{o}_b(V) & \forall n\in N\\
& & &\mathbf{V_n^{min}}\leq |V_n|\leq \mathbf{V_n^{max}} & \forall n\in N\\
& & & \mathbf{P_g^{min}}\leq Re(S_g^{gen})\leq \mathbf{P_g^{max}} & \forall g\in G\\
& & & \mathbf{Q_g^{min}}\leq Im(S_g^{gen})\leq \mathbf{Q_g^{max}} & \forall g\in G\\
& & & |S^{o}_b(V)|\leq \mathbf{S_l^{max}} & \forall b\in B\\
& & & |S^{d}_b(V)|\leq \mathbf{S_l^{max}} & \forall b\in B
\end{aligned}
\end{equation} where the constants are in bold. $Re(z)$ stands for the real part of $z\in \mathbb{C}$ and $Im(z)$ for the imaginary part. $S^{load}_n$ is a constant and represents the electrical demand (or load) at bus $n$. $B^-(n)$ is the set of entering branches at bus $n$ and $B^+(n)$ the set of exiting branches. Each function $c_g$ is quadratic. The power $S^{shunt}_n$ is defined as $S^{shunt}_n=-b_n |V_n|^2$. The powers at the origin and at the destination of a branch $b\in B$ are defined from the voltages as following:
\begin{equation}
S^{o}_b(V)=V_{o(b)}(\overline{(Y_b)_{11} V_{o(b)}+ (Y_b)_{12}V_{d(b)}})
\label{Sorig_V}
\end{equation}

\begin{equation}
S^{d}_b(V)=V_{d(b)}(\overline{(Y_b)_{21} V_{o(b)}+ (Y_b)_{22}V_{d(b)}})
\label{Sdest_V}
\end{equation} where $o(b)$ stands for the origin of branch $b$, $d(b)$ for the destination of branch $b$ and $Y_b$ is the admittance matrix of the branch $b\in B$ defined from physical characteristics such as the susceptance (see \cite{HistoryofOPF} for more details).\\

This mathematical formulation contains two types of complex variables:

\begin{itemize}
	\item a voltage $V_n\in \mathbb{C}$ for each bus $n\in N$   
	\item a power $S^{gen}_g\in \mathbb{C}$ for each generator bus $g\in G \subset N$
\end{itemize}
The objective is to minimize the real generation costs subject to power balance constraints at every bus plus several safety constraints: the voltage magnitude is bounded at each bus, active and reactive powers at each generator bus are also bounded and finally there are thermal limits for each branch. The next section demonstrates how ACOPF problems can be constructed with our module.

\subsection{Modelling flexibility}

To meet the need for modelling flexibility, we have designed a generic structure where the user can define network elements himself. In this structure, the network is composed of a set of buses and a set of links between buses. Then several elements can be associated to each bus and to each link. For example, a shunt element can be associated to a bus and a $\Pi$ transmission line can be associated to a link. One of the advantages of this structure is to allow the association of several elements to the same bus or to the same link. For example, one bus can be connected to a shunt element, to a load element but also to several generator elements and one link can be associated to two different $\Pi$ transmission lines. The user can define a bus element of the power network by specifying:
\begin{itemize}
\item the variables this element needs
\item its contribution into the power balance
\item the constraints associated with this element
\item its contribution into the cost function
\end{itemize}

This way of defining the elements enables an automatical writing of power balance constraints, which basically consists in going through all the elements associated to a node and summing up their contribution. 

For example, to define \eqref{ACOPF}, four bus elements are needed. These bus elements are summarized in Table \ref{MATPOWERbuselements}. The first one is the voltage element and it defines a voltage variable and bound constraints on the voltage magnitude. There is no contribution into the power balance or the cost function for a voltage element. The next two are the load and the shunt element. Both only have a contribution in the power balance as defined in \eqref{ACOPF}. Finally, a generator element defines a power variable along with bounds constraints on its real and imaginary part. The contribution into the power balance constraint is the opposite of the power following the model \eqref{ACOPF}. The generator element is the only element having a contribution in the cost function: there is a quadratic cost depending on the real part of the power variable.

Similarly, a link element is defined by specifying:
\begin{itemize}
\item the variables this element needs
\item the power at the origin of the link
\item the power at the destination of the link
\item the constraints associated with this element
\item its contribution into the cost function
\end{itemize}

To define \eqref{ACOPF}, only one link element is needed: the $\Pi$ transmission line described in Table \ref{MATPOWERlinkelements}. This element involves two voltage variables: one for the origin bus and one for the destination bus. The powers at the origin and at the destination are respectively defined in \eqref{Sorig_V} and \eqref{Sdest_V}. Finally, it involves thermal limits constraints.

\begin{table}[!ht]
\renewcommand{\arraystretch}{1.2}
\caption{Bus elements}
\label{MATPOWERbuselements}
\noindent
\centering
    \begin{minipage}{\linewidth} 
    \renewcommand\footnoterule{\vspace*{-5pt}} 
\begin{center}
\resizebox{\textwidth}{!}{
\begin{tabular}{|c|c|c|c|c|}
  \hline
  Element & Variables  & Power balance & Constraints & Cost\\
  \hline
   Voltage & $V_n$  & - & $\mathbf{V_n^{min}}\leq |V_n|\leq \mathbf{V_n^{max}}$ & - \\
  \hline
  Load & -  & $\mathbf{S_n^{load}}$ & - & - \\
  \hline
   Shunt & -  & $-\mathbf{b_n} |V_n|^2$ & - & - \\
  \hline
   Generator & $S^{gen}_n$ & $-S^{gen}_n$  & $\mathbf{P_g^{min}}\leq Re(S_g^{gen})\leq \mathbf{P_g^{max}}$&  $c_n(Re(S_n^{gen})$\\
& & & $\mathbf{Q_g^{min}}\leq Im(S_g^{gen})\leq \mathbf{Q_g^{max}}$ & \\
  \hline
\end{tabular}}
\end{center}
    \end{minipage}
\end{table}

\begin{table}[!ht]
\renewcommand{\arraystretch}{1.2}
\caption{Link elements}
\label{MATPOWERlinkelements}
\noindent
\centering
    \begin{minipage}{\linewidth} 
    \renewcommand\footnoterule{\vspace*{-5pt}} 
\begin{center}
\resizebox{\textwidth}{!}{
\begin{tabular}{|c|c|c|c|c|c|}
  \hline
  Element & Variables  & Power at origin & Power at destination & Constraints & Cost\\
  \hline
  $\Pi$ transmission line & $V_{o(b)}$,$V_{d(b)}$  & $S^{o}_b(V)$& $S^{d}_b(V)$ & $|S^{o}_b(V)|\leq \mathbf{S_l^{max}}$ &-\\
   & & & & $|S^{d}_b(V)|\leq \mathbf{S_l^{max}}$ &\\
  \hline
\end{tabular}}
\end{center}
    \end{minipage}
\end{table}

We have tested our module on several MATPOWER instances (without thermal limits constraints) to validate our approach. The results are not detailed for brevity's sake but the next section presents the more general PSCOPF problem for which computational results are given in section D. 
 
\subsection{Preventive-Security Constrained Optimal Power Flow}
The ACOPF problem only focuses on the base case network but as grid flexibility, reliability and safety are taking greater importance within the context of energy transition, grid failure is anticipated in more general problems such as the PSCOPF problem. Besides inheriting the OPF’s nonconvexity, this problem anticipates credible contingency cases, which means it considers simultaneously several configurations of the power network. A contingency case is defined by the loss of a network element: a generator, a transmission line or a transformer. The PSCOPF involves several nonseparable variants of the OPF problem in the same optimization problem. In other words, it is not possible to solve the ACOPF problems separately as all contingency dispatches depend on the base case optimal dispatch. Moreover, this dependance involves the use of binary variables in the model.\\

The PSCOPF model given in the Beta Phase of the Grid Optimization Competition can be written as follows:
\begin{equation}\label{PSCOPF}\tag{PSCOPF}
\begin{aligned}
&\underset{V_k, S^{gen}_k, \Delta_k}{\text{min}}
\sum_{g\in G} c_g(Re(S_{g,0}^{gen}))\\
&S^{gen}_{n,k} = \mathbf{S^{load}_n}+S^{shunt}_n+\\
&\sum_{b\in B^-(n)}S^{d}_b(V_k) +\sum_{b\in B^+(n)}S^{o}_b(V_k)& \forall n\in N, \forall k\in K\\
&\mathbf{V_n^{min}}\leq |V_{n,k}|\leq \mathbf{V_n^{max}} & \forall n\in N, \forall k\in K\\
& \mathbf{P_g^{min}}\leq Re(S_{g,k}^{gen})\leq \mathbf{P_g^{max}} & \forall g\in G, \forall k\in K\\
& \mathbf{Q_g^{min}}\leq Im(S_{g,k}^{gen})\leq \mathbf{Q_g^{max}} & \forall g\in G, \forall k\in K\\
& |S^{o}_b(V_k)|\leq \mathbf{S_l^{max}} & \forall b\in B, \forall k\in K\\
& |S^{d}_b(V_k)|\leq \mathbf{S_l^{max}} & \forall b\in B, \forall k\in K\\
& Re(S_{g,k}^{gen})=Re(S_{g,0}^{gen})+\bm{\alpha_g} \Delta_k &\forall g\in G, \forall k\in K^*\\
&|V_{g,k}| <|V_{g,0}| \Rightarrow Im(S_{g,k}^{gen})=\mathbf{Q_g^{max}}& \forall g\in G, \forall k\in K^*\\
&|V_{g,k}| > |V_{g,0}| \Rightarrow Im(S_{g,k}^{gen})=\mathbf{Q_g^{min}} &\forall g\in G, \forall k\in K^*
\end{aligned}
\end{equation} where $K=\{0,1,2,\ldots, N_K\}$ represents the different network cases: the index $0$ represents the base case while the other indices represent contingency cases. The constant $\alpha_g$ is the participation share of generator $g$ and the variable $\Delta_k$ is the pre-recovery real power shortfall in case $k$. The other notations are all defined in section A. 

The last two constraints are called PV/PQ switching constraints: the generator reactive power adjusts to maintain voltage magnitude from the base case. In other words, a decrease in voltage magnitude is compensated by a maximal injection of reactive power and vice-versa. We model these constraints using binary variables.

PSCOPF problems can be constructed with our module thanks to modelling flexibility. Computational results for Grid Optimization Competition datasets are presented in the next section.

\subsection{Computational results}

We have tested our module on the datasets given in the Beta Phase of the Grid Optimization Competition \cite{GOC}. While the Beta Phase is over, several PSCOPF datasets for small networks were made available during the Beta Phase along with their solutions. Some are still available at the address \cite{GOC_datasets_Beta}. These datasets are briefly described in Table \ref{tabdatasets}. There are four datasets which all represent small networks: IEEE14 and Modified\_IEEE14 datasets represent 14-buses networks and RTS96 and Modified\_RTS96 represent 73-buses networks. The main difference between a dataset and its modified version is the number of branches: the modified version is the same network but with one fewer branch. For IEEE14 and Modified\_IEEE14 datasets, there is only one contingency to consider, which corresponds to the loss of a transmission line. For RTS96 and Modified\_RTS96 datasets, there are respectively ten and nine contingencies: one out of ten (or nine) branches can be lost. Finally, in each dataset, there are 100 scenarios, i.e., 100 different instances. The differences can be related to electrical demand for example.

\begin{table}[!ht]
\renewcommand{\arraystretch}{1.2}
\caption{Description of datasets}
\label{tabdatasets}
\noindent
\centering
    \begin{minipage}{\linewidth} 
    \renewcommand\footnoterule{\vspace*{-5pt}} 
\begin{center}
\resizebox{\textwidth}{!}{
\begin{tabular}{|c|c|c|c|c|}
  \hline
  Network & $\#$ nodes  & $\#$ lines & $\#$ contingencies & $\#$ scenarios\\
  \hline
  IEEE14 & 14  & 20 & 1 & 100 \\
  \hline
   Modified\_IEEE14 & 14  & 19 & 1 & 100 \\
  \hline
  RTS96 & 73  & 120  & 10 & 100 \\
  \hline
   Modified\_RTS96 & 73  & 119  & 9 & 100 \\
  \hline
\end{tabular}}
\end{center}
    \end{minipage}
\end{table}

We construct the PSCOPF problems as $(MIPOP-\mathbb{C})$ using our Julia package. They are then converted to Mixed-Integer NonLinear Problems (MINLP) using the rectangular form. As the problem combines difficulties, we apply a three-step method to find a good feasible solution. In the first step, the continuous relaxation of the problem is solved, which amounts to finding solutions for all OPF problems as if they were independent. The goal of the second and third steps is to reach convergence with binary variables. In the second step, the MINLP problem is solved using the Knitro \cite{KNITRO} solver complementarity option, that is, the problem solved is the continuous reformulation obtained by converting all binary variables into complementary constraints and adding them to the objective function with a penalization parameter. This parameter is automatically updated by Knitro. Finally, the binary variables are fixed in the third step using the second step solution, which ensures to get a feasible solution. 

The tests were carried out on a Processor Intel® Core™ i7-6820HQ CPU @2.70GHz using our Julia module with Julia Version 0.6.1, AMPL [Version 20161231] and the solver ArtelysKnitro\cite{KNITRO} 10.3.0. The feasiblity tolerance  (feastol) was $10^{-6}$ and the optimality tolerance (optol) was $10^{-3}$. 

We say that an instance is solved if we can compute a feasible solution as good as the one given in the datasets with respect to optimal value and feasibility. The results are summarized in Table \ref{tabresults} and were confirmed by the test platform of the Grid Optimization Competition. Numerical issues have arised in this work. Indeed, a large number of instances are solved but not all of them: 90/100 for IEEE14 and RTS96, 89/100 for Modified\_RTS96 and only 84/100 for Modified\_IEEE14. The continuous relaxation (step 1) always converges but reaching convergence with the binary variables is heavily dependent on scaling: with different scalings but the same stopping criteria, the number of solved instances differ and in addition the instances which are not solved are not always the same. For example, for IEEE14, only 78 out of 100 instances are solved without scaling.
 
\begin{table}
\renewcommand{\arraystretch}{1.3}
\caption{Results}
\label{tabresults}
\noindent
\centering
    \begin{minipage}{\linewidth} 
    \renewcommand\footnoterule{\vspace*{-5pt}} 
\begin{center}
\resizebox{\textwidth}{!}{
\begin{tabular}{|c|c|c|c|c|c|c|c|}
  \hline
  Network & $\#$variables  & $\#$constraints & $\#$nonzeros & $\#$nonzeros & $\#$solved \\
   & & & Jacobian & Hessian & scenarios\\
  \hline
  IEEE14 & 92  & 207 & 937 & 245 & 90/100\\
  \hline
  Modified\_IEEE14 & 92  & 203  & 905 & 237 & 84/100\\
  \hline
  RTS96 & 4784  & 12157 & 49838 & 7199 &  90/100\\
  \hline
   Modified\_RTS96 &  4340 & 10987 & 44960 & 6512 & 89/100\\
  \hline
\end{tabular}}
\end{center}
    \end{minipage}
\end{table}

The application of MathProgComplex.jl to ACOPF and PSCOPF problems has validated our module while highlighting future research directions that are presented in the next section.

\section{Future research}

There is still a lot of research to be done for $(POP-\mathbb{C})$. Three topics are to be explored to further improve the capabilities of our Julia module: the computation of lower bounds, the numerical conditioning and the design of algorithms for problems with complex variables.

Studies carried out by Josz and Molzahn made the computation of lower bounds for $(POP-\mathbb{C})$ a promising research direction \cite{josz2018lasserre}. The computation of these bounds is indeed essential for $(POP-\mathbb{C})$ since most of these problems are nonconvex, which means that only local solutions can be provided by solvers. To evaluate the quality of local solutions, lower bounds have to be computed by solving convex relaxations. For instance, local solutions computed for ACOPF problems are often optimal or near-optimal but some convex relaxations have to be solved to prove it \cite{Zero_duality_gap,lu2018tight,coffrin2017,Bingane_conic_relax}. The Lasserre hierarchy provides a convergent sequence of Semi Definite Programming (SDP)\cite{SDP} relaxations for polynomial optimization in complex or real variables \cite{Josz_complex_SOS}. As a first tool to compute lower bounds, a Lasserre hierarchy implementation is integrated into MathProgComplex.jl, for complex or real problems. A few options are available: sparsity of problems can be exploited (the set of exponents can be split into smaller cliques) and it is also possible to apply different orders on the different constraints. Finally, some symmetries can be speficied to simplify the problems. As alternative, other conic relaxations such as Second Order Cone Programming (SOCP) ones \cite{SOCP_relax_Sun_2017} could be implemented in the future. 

The numerical conditioning for $(POP-\mathbb{C})$ is another issue to tackle since conditioning nonconvex problems is indispensable, as shown in section III.D. The question then arises: can we extend to complex problems the work done in the real case \cite{conditioning}?

Currently, MathProgComplex.jl has only been tested on ACOPF and PSCOPF problems but	other classes of power system problems could be implemented, like Unit Commitment problems. More broadly speaking, our module is a tool for any $(POP-\mathbb{C})$ regardless of the field. The ultimate purpose is to design algorithms for problems with complex variables. Even if there are not many papers on the subject, the researchers familiar with it seem to encourage the use of the complex structure. For instance, Gilbert and Josz push towards a solver for complex SDP problems in \cite{jcg_SDP_solver_complex_numbers} and Chen, Atamt{\"u}rk and Oren take profit of the complex structure to get a tighter relaxation in \cite{B&C_for_nonconvex_QCQP_varCbounded}.

\section{Conclusion}

As accurate resistive losses computation is only possible with Alternating Current (AC) modelling, most power system problems involve complex variables. To overcome the lack of tools for problems with complex variables, we proposed a Julia module for $(MIPOP-\mathbb{C})$. This module was validated on ACOPF and PSCOPF problems. To continue the work on power systems, this module could be integrated into the module PowerModels.jl to provide a complex formulation of the ACOPF problem.





%

\end{document}